\documentclass[12pt,leqno]{article} 

\usepackage{amsmath} 
\usepackage{amscd} 
\usepackage{amssymb} 
\usepackage{amsfonts}
\usepackage{amsthm}
\usepackage{makeidx} 
\usepackage{verbatim}

\newlength{\fixboxwidth}
\newcommand{\tabmath}[1]{%
	\setlength{\fboxrule}{0pt}%
	\fbox{${#1}$}%
}

\newcommand{\R}{{\mathbb R}}
\newcommand{\N}{{\mathbb N}}
\newcommand{\SP}{{\mathbb S}}
\newcommand{\PP}{{\mathbb P}}

\newcommand{\bi}{{\bf{i}}}
\newcommand{\bx}{{\bf{x}}}

\newcommand{\refx}{\medskip\noindent}
\renewcommand{\rho}{{\varrho}}
\def\min{{\rm min}}  
\def\old{{\rm old}} 
 
\def\a{{\alpha }} 
\def\g{{\gamma }} 

\newcommand{\wt}{\widetilde }

\theoremstyle{plain}
\newtheorem{theorem}{Theorem}

\newtheorem{lemma}{Lemma}

\theoremstyle{definition}
\newtheorem{rem}{Remark}

\begin{document}

\title{Cubature Formulas for Symmetric Measures in Higher Dimensions
with Few Points}

\author{Aicke Hinrichs\footnote{Research of the first author 
was supported by the DFG Emmy-Noether grant 
Hi 584/2-4.}, Erich Novak \\[5pt] 
\normalsize  Mathematisches Institut, Universit\"at Jena, \\
\normalsize   Ernst-Abbe-Platz 2, D-07743 Jena, Germany. \\
\normalsize   Email: hinrichs@math.uni-jena.de, novak@math.uni-jena.de} 

\maketitle

\begin{abstract}
We study cubature formulas for $d$-dimensional integrals with an 
arbitrary symmetric weight function of 
product form. We present
a construction that yields a high polynomial exactness:
for fixed degree $\ell=5$ or $\ell=7$
and large dimension $d$ the number of knots
is only slightly larger than the lower bound of M\"oller 
and much smaller compared to the known constructions. 

We also show, for any odd degree $\ell = 2k+1$, that the minimal 
number of points is almost independent of the weight function. 
This is also true for the integration over the (Euclidean) sphere.
\end{abstract}

\maketitle

{\small 


 } 

\section{Introduction}\label{s1} 

Let us start with a special case of our results:
We find cubature formulas
%
%
with 
$$
N(5, d,  1 ) = d^2 + 7d+1 ,  \quad \hbox{and} \quad 
N(7, d,  1 ) = (d^3+ 21d^2 +20d+3 )/3
$$
points such that 
the integral 
$$
I_d (f) = \int_{[-1,1]^d} f( \bx ) \,  d\bx 
$$
is exactly computed for all polynomials of degree 
at most 5 or 7, respectively. 
This improves the known cubature formulas 
for degree $5$ and $d\ge 8$ and for degree $7$ with 
$d \ge 10$. 
The lower bound of M\"oller (1979)  takes the form
\begin{equation}   \label{moeller}  
N_\min  (5, d, 1 ) \ge d^2  +d +1 \qquad \hbox{and} \qquad
N_\min  (7, d, 1 ) \ge (d^3+ 3d^2 + 8d )/3 . 
\end{equation} 
Hence, for our method, we obtain 
\begin{equation}  
N (5, d, 1 ) \approx N_\min (5,d, 1 )
\qquad \hbox{and} \qquad
N (7, d, 1 ) \approx N_\min (7,d, 1 )  . 
\end{equation}  
We use $\approx$ to denote the strong 
equivalence of sequences, i.e., 
\[
v_n \approx w_n \qquad \text{iff} \qquad \lim_{n\to \infty }v_n/w_n=1.
\]
The best results (for large $d$) from the literature, 
see   Stroud (1971)
and the online tables of Cools, 
see Cools (2003), 
are given by 
\begin{equation}   \label{old}  
N_\old  (5, d, 1 ) = 2 d^2 +1 \qquad \hbox{and} \qquad
N_\old  (7, d, 1 ) = (4d^3 - 6d^2 +14 d +3)/3 . 
\end{equation}  

More generally, we study cubature formulas 
\begin{equation}\label{g1}
Q_n(f) = \sum_{i=1}^n a_i\, f(\bx_i),  
\qquad a_i \in \R, \ \bx_i \in \Omega ,
\end{equation}
for $d$-dimensional integrals
\begin{equation}\label{g2}
I_d^\rho  (f) = \int_{\Omega} f( \bx ) \,  \rho  (\bx )  \, d\bx .
\end{equation}
Concerning the integral we always assume
$$
\Omega = \Omega_1 \times \dots \times \Omega_d
$$
with symmetric (and possibly unbounded) intervals $\Omega_j \subset \R$
and the  
product form
$$
\rho  (\bx ) = \rho_1 (x_1) \dots \rho_d (x_d) 
$$
of the weight function $\rho$. 
We assume that the $\rho_i$ are symmetric, 
$$
\rho_i (x) = \rho_i (-x)
$$
with $\rho_i \geq 0$ and integrability of all polynomials, 
although these assumptions can be relaxed.
Some of our results can be slightly improved in the fully symmetric 
case where, in addition, all the $\rho_i$ coincide. 

Let $\PP(\ell,d)$ be the space of all polynomials in $d$ variables 
of (total) degree at most~$\ell$. A cubature formula $Q_n$ has a
degree $\ell$ of exactness if
$$
\phantom{\qquad \forall \, f \in \PP(\ell,d).}
Q_n(f) = I_d^\rho (f), 
\qquad \forall \, f \in \PP(\ell,d).
$$
We define
$$
N_{\rm min} (\ell , d, \rho )
$$
to be the minimal number $n$ of knots needed by any cubature
formula $Q_n$ of degree $\ell$ of exactness. 

The numbers $N_{\rm min}(\ell,d, \rho)$ and corresponding cubature
formulas are only known in exceptional cases, see, e.g., Schmid
(1983), Berens, Schmid, Xu (1995), and Cools (1997).
Thus one is interested in upper and lower bounds 
for this quantity. 

One is often interested in cubature formulas
with knots inside the domain and positive weights. While
$\bx_i \in \Omega$ can always be satisfied by our method, we usually
have positive and negative weights. 
Actually we request $\bx_i \in \Omega$, see \eqref{g1}, although 
the lower of M\"oller also holds without this assumption. 

\section{Problem, Main Results, and Conjecture} 

The lower bound of M\"oller (1979) 
for centrally symmetric weight functions is 
the following: 
If $k$ is odd then 
$$
N_\min (2k+1, d, \rho ) \ge
2 \dim \PP_e (k,d) = \binom{d+k}{d} +  \sum_{s=1}^{d-1} 
2^{s-d} \binom{s+k}{s} .
$$
If $k$ is even then 
$$
N_\min (2k+1, d, \rho ) \ge
2 \dim \PP_o (k,d) -1 = \binom{d+k}{d} +  \sum_{s=1}^{d-1} 
(1- 2^{s-d})   \binom{s+k-1}{s} .
$$
Here $\PP_e(k,d)$ denotes the subspace of $\PP(k,d)$ generated
by even polynomials and $\PP_o(k,d)$ is the subspace generated 
by odd polynomials. 
We obtain \eqref{moeller} as special cases
and for large $d$ the lower bounds are of the order 
$$
\approx  \frac{ 2 d^k}{k!} .
$$
See the book 
Mysovskikh (1981) or Cools (1997) and, 
for the explicit formula, Lu, Darmofal (2004). 

The best upper bounds were of the form 
\begin{equation}   \label{eq15} 
\approx  \frac{2^k \, d^k}{k!} .
\end{equation} 
They can be proved with 
``fully symmetric formulas'' (if the $\rho_i$ are equal) 
or (in the general case) with the ``Smolyak method'' 
or with ``sparse grids''. All these notions are very much 
related, see Section~\ref{smo}.  
Even for special weight functions $\rho$ and/or 
for special $\ell=2k+1$ better bounds were not known. 
Hence there is a gap between the lower and the upper bound 
of a factor of $2^{k-1}$
and we only knew (before we wrote this paper) of one exception: 
For the weight function 
\begin{equation}   \label{eq16} 
\rho (\bx ) = \exp({-\Vert \bx \Vert_2^2}) ,
\end{equation} 
it is known for  $\ell =5$ that
\begin{equation}   \label{eq17} 
d^2+3d+3
\end{equation} 
function values are enough, see Lu, Darmofal (2004). 

Observe that the weight function \eqref{eq16} is 
invariant with respect to rotations. 
Hence one might ask whether a result similar to 
\eqref{eq17} 
holds for all symmetric weight functions. 
We conjecture that 
\begin{equation}   \label{eq18}  
N_\min (2k+1,d ,\rho) \approx \frac{ 2 d^k}{k!} 
\end{equation} 
holds for all $\rho$ and all $k$, hence the M\"oller bound 
is almost optimal. 
In this paper we prove this conjecture for $k=2$ and $k=3$, 
see Theorem 1 for more details. 
We also prove that the numbers $N_\min (2k+1, d, \rho)$ 
only mildly depend on the weight function $\rho$,
see Theorem 2 for the details. 

\section{Some facts about the Smolyak method} \label{smo} 

We study a special case of the Smolyak method, 
as we need it in the following. 
We also present methods with the upper bound \eqref{eq15},
since they are used (twice) for our new 
algorithm with the improved bound. 
We believe that this proof technique can be used to 
establish the conjecture \eqref{eq18} in full generality. 
Everything in this section is known or a minor modification of known 
results, see Novak, Ritter (1999). 

We construct  cubature formulas to compute the
integral \eqref{g2} as follows.  
First we select 
quadrature formulas 
$U^1_j, U^2_j, \dots $ to compute the one-dimen\-sional integrals
\[
\int_{\Omega_j} f(x) \,  \rho_j (x) \, dx.
\]
These formulas should have the following properties:
The formula $U^i_j$ is exact for all univariate polynomials of 
degree $m_i$, where 
\begin{equation}  \label{eq19}   
m_i \ge  2i-1 .     
\end{equation} 
The formula $U^i_j$ uses the knots $X^i_j$, 
the number $n_i = | X^i_j| $ of knots satisfies
\begin{equation} \label{eq20}  
n_i \le 2i-1 .
\end{equation} 
We also assume that the $X^i_j$ are symmetric and ``embedded'' or ``nested'', 
i.e., 
\begin{equation} \label{eq20a}   
X^{i-1}_j \subset X^i_j  \qquad \hbox{  for every $i$ and $j$} .
\end{equation}  

By \eqref{eq19} and  \eqref{eq20} the weights of $U^i_j$ 
are uniquely determined by 
its knots. Formulas with this property are often called 
interpolatory quadrature formulas. 
For simplicity we assume in this paper 
that the numbers $m_i$ and $n_i$ do not depend on the coordinate $j$. 
The formula $U^i_j$, however, may depend on $j$. 

A product formula $U^{i_1}_1 \otimes \dots \otimes U^{i_d}_{d}$ 
needs $n_{i_1} \dots n_{i_d}$ 
function values, sampled on a grid. The Smolyak formulas
$A(q,d)$ are linear combinations of product formulas with
the following key properties. Only products with a relatively 
small number of knots are used and the linear combination is chosen in 
such a way that the interpolation property for $d=1$ is preserved 
for $d>1$. The formula $A(q,d)$ is defined by 
\begin{equation}  \label{Aqd} 
A(q,d) = 
\sum_{q-d+1 \le |\bi| \le q} (-1)^{q-|\bi|} \cdot 
{d-1  \choose  q-|\bi| } \cdot
(U^{i_1}_1 \otimes \cdots \otimes U^{i_d}_d), 
\end{equation}  
where $q \geq d$, $\bi \in \N^d$, and $|\bi|=i_1 + \dots + i_d$. 

The cubature formula $A(q,d)$ is based on the sparse grid
$$
H(q,d) =
\bigcup_{|\bi| = q} 
X_1^{i_1} \times \cdots \times X_d^{i_d}, 
$$
we use
\[
n = n(q,d)
\]
to denote the cardinality of $H(q,d)$.\footnote{Observe 
that some elements of the sparse grid might get a zero weight 
in the formula $A(q,d)$. This would decrease the number
of needed function values. 
Hence the ``actual'' number of needed function values 
for $A(q,d)$ might be smaller than $n(q,d)$.}
In particular we have $n(q,1) = n_q$ and we put $n(0,1)=n_0=0$.
The recursion formula 
\begin{equation}\label{g30}
n(q+1,d+1) = \sum_{s=1}^{q-d+1} n(q+1-s, d) \cdot
(n_s-n_{s-1})
\end{equation}
for $n(q,d)$ is known, see Novak, Ritter (1999). 

\begin{rem} 
Cubature formulas with high polynomial exactness 
are not often used if $d$ is large, say $d>5$. 
One major exception is the class of fully symmetric rules for  
the fully symmetric case, where also 
$$
\rho_1 = \dots = \rho_d . 
$$
Fully symmetric cubature formulas were
developed by
Lyness (1965a, 1965b), 
McNamee and Stenger (1967), Genz (1986), Cools and 
Haegemans (1994), Capstick and Keister (1996),
Genz and Keister (1996) and other authors. 
The best results with respect to polynomial exactness are obtained by 
Genz (1986) and Genz and Keister (1996). 
The fully symmetric formulas from 
Genz (1986) and Genz and Keister (1996)
are of the Smolyak form \eqref{Aqd}.
Numerical integration with the Smolyak construction 
was already studied in Smolyak (1963).
There are many other papers on the Smolyak method. 
The papers 
Gerstner, Griebel (1998),
Novak, Ritter (1999),
and Petras (2003) 
study the polynomial exactness of $A(q,d)$. 
See also 
Novak, Ritter, Schmitt, Steinbauer (1999) and 
the recent survey on sparse grids by 
Bungartz, Griebel (2004).
\end{rem} 

The following result is well known, see Corollary 1 of 
Novak, Ritter (1999). 

\begin{lemma}  \label{lemma1}  
Assume \eqref{eq19}. 
Then $A(d+k,d)$ has 
(at least) a degree $\ell=2k+1$ of exactness.
\end{lemma}

Now we present formulas for the number $n(q,d)$ of knots that 
are used by $A(q,d)$. We consider two cases, important for 
the following. 

\subsubsection*{The case $n_i=2i-1$.} 

Using \eqref{g30} one obtains the recursion
\begin{equation}         \label{zahl1}  
n(q+1,d+1) = n(q,d+1) + n(q,d) + n(q-1,d)
\end{equation}  
for $q \geq d$ and $n(q,1)= 2q-1$ and $n(d,d)=1$.  
Table~1 consists of numbers $n(q,d)$ with minimal $q$ such
that $n(q,d) \geq \ell$, these numbers are called 
$N(\ell,d)$.

\begin{table}[h]
\caption{Number of knots for Smolyak's method with $n_i=2i-1$}
\label{ta3}
\[
\begin {array}{rrrrrr} 
\ell  & N(\ell ,5) & N (\ell ,10) & N (\ell ,15) & N (\ell ,20) 
& N (\ell ,25)  \\[\medskipamount]
3&11&21&31&41&51
\\\noalign{\medskip}
5&61&221&481&841&1\,301
\\\noalign{\medskip}
7&231&1\,561&4\,991&11\,521&22\,151
\\\noalign{\medskip}
9&681&8\,361&39\,041&118\,721&283\,401
\\\noalign{\medskip}
11&1\,683&36\,365&246\,047&982\,729&2\,908\,411
\\\noalign{\medskip}
13&3\,653&134\,245&1\,303\,777&6\,814\,249&24\,957\,661
\\\noalign{\medskip}
15&7\,183&433\,905&5\,984\,767&40\,754\,369&184\,327\,311
\\\noalign{\medskip}
17&13\,073&1\,256\,465&24\,331\,777&214\,828\,609&1\,196\,924\,561
\end {array}
\]
\end{table}

Using \eqref{zahl1} one can get an explicit formula 
for $n(k+d,d)$, see Novak, Ritter (1999). 


\begin{lemma}\label{t2}
For every $k \in \N_0$ and $d \in \N$ we have
$$
n(k+d,d) = \sum_{s=0}^{\min(k,d)} \binom{k}{s}
\cdot \binom{k+d-s}{k}. 
$$
\end{lemma}

\begin{rem}\label{r5}
Lemma~\ref{t2} immediately implies
\begin{equation}\label{g24}
n(k+d,d) \leq \binom{k+d}{d} \cdot \min (2^k,2^d).
\end{equation}
\end{rem}

\subsubsection*{The case $n_i=2i-1$ for $i \not= 3$ and $n_3=3$.} 


If we take the Gaussian formulas $U^2_j$ with 3 knots for
$\varrho_j$, then we already have exactness 5 and so we can take 
$U_j^3=U^2_j$ and still have \eqref{eq19}. Altogether 
we have
\begin{equation}  \label{ass2} 
n_i=2i-1 \quad \hbox{for} \quad  i \not= 3, \  \   n_3=3. 
\end{equation} 
Observe that in this case the sets $X^2_j$ are determined by the weights 
$\varrho_j$, we cannot choose these sets. 
All the other sets $X^i_j$ can be chosen arbitrarily for 
$i >2$, but we still assume \eqref{eq20a}. 
Similarly as \eqref{zahl1} we now obtain from 
\eqref{g30}  the recursion
\begin{eqnarray*} 
 n(q+2,d+1) &=& n(q+1,d+1)+n(q+1,d)+n(q,d) 
  \\ \nonumber && -2n(q-1,d)+4n(q-2,d)-2n(q-3,d).
\end{eqnarray*}  
With this simple modification we obtain the values of Table~2.

\begin{table}[h]   
\caption{Number of knots for method \eqref{ass2}}
\label{ta4}
\[
\begin{array}{rrrrrr} 
\ell  & N(\ell ,5) & N (\ell ,10) & N (\ell ,15) & N (\ell ,20) 
& N (\ell ,25)  \\[\medskipamount]
3&11&21&31&41&51
\\\noalign{\medskip}
5&51&201&451&801&1\,251
\\\noalign{\medskip}
7&151&1\,201&4\,151&10\,001&19\,751
\\\noalign{\medskip}
9&401&5\,301&27\,701&90\,601&227\,001
\\\noalign{\medskip}
11&1\,003&19\,505&146\,507&643\,009&2\,040\,011
\\\noalign{\medskip}
13&2\,133&63\,805&655\,017&3\,775\,769&15\,056\,061
\\\noalign{\medskip}
15&4\,223&188\,745&2\,584\,167&19\,111\,089&94\,680\,111
\\\noalign{\medskip}
17&8\,113&511\,625&9\,224\,937&85\,920\,449&522\,028\,561
\end {array}
\]
\end{table}

\begin{rem} 
Later the following will be important for the two versions 
of Smolyak's algorithm: 
In the case $n_i=2i-1$ we can take arbitrary symmetric sets 
$X^i_j$, in particular we can take
$$
X^2_1 = \dots = X^2_d .
$$
We also can normalize the weights $\varrho_j$ in such a way 
that the $U^2_j$ have the form
$$
U^2_j (f) = \gamma f(-x) + \beta_j f(0) + \gamma f(x) ,
$$
where $\gamma$ (and $x$) do not depend on $j$. 
In addition, 
we can choose the $X^i_j$ in such a way that 
$\Vert \bx \Vert_2 \le \alpha$ for
each $\bx \in  H(q,d)$, where  
$\a$ is the (given) radius 
of the domain $\Omega$ of integration. 
This means that each rotation maps $\bx$ to a point in 
$\Omega$. 

In the second case, however, we have to use 
the 3 Gau\ss -knots for $X^2_j=X^3_j$. 
\end{rem} 

\begin{rem}    \label{rem4}
Later we project the points $H(q,d)$ of $A(q,d)$ to 
a sphere of fixed radius.
The origin is not projected. 
This projection reduces the number 
of points, the number of projected points 
$n^*(d+k,d)$  also depends 
on the sets $X^i_j$.
We only need the second case, where $n_i = 2i-1$ for $i\neq 3$ and $n_3=3$.
In the case $k=2$ and $k=3$ one obtains 
$$
n(d+2,d)= 2d^2+1 \quad \hbox{and} \quad n^*(d+2,d)= 2d^2
$$
and
$$
n(d+3,d)= (4d^3-6d^2+20d+3)/3  \quad \hbox{and} \quad n^*(d+3,d)
= (4d^3-6d^2+8d)/3 . 
$$
For the last formula observe that $H(d+3,d)$ contains 7 points 
of the form $\bx = ( \alpha, 0, \dots , 0)$ that are projected 
onto two different points, hence 
$$
n(d+3,d) = n^*(d+3,d) + 4d+ 1 .
$$
It seems to be difficult to compute the smallest possible 
number $n^*(d+k,d)$ for general $k$, but it is clear that 
$$
n(d+k,d) \ge n^*(d+k,d) \ge 2^k \binom{d}{k}.
$$
Hence, for large $d$, we have
$n(d+k,d) \approx n^*(d+k,d)$.
\end{rem} 

\section{Known results for the Lebesgue measure} 

Here we explain the best known upper bounds for 
$N_\min (\ell, d, 1)$ that we found in the literature. 
Again we only discuss results for large $d$.\footnote{We illustrate 
this by an example. In the case $d=10$ and $\ell=13$ we will 
mention a method of Genz (1986) using $n=60\,205$ function values. 
In the same paper Genz presents another method using only
$n=37\,389$. This method, however, uses more than $2^d$ points 
for general $d$ and hence is not good for ``large'' $d$.} 

The results for $\ell \in \{ 3, 5, 7 \}$ are classical results 
that can be found in Stroud (1971): 

$n=2d$ \quad  for the degree $\ell=3$;
this bound is sharp, $N_\min (3,d,\rho) = 2d$; 

$n=2d^2+1$ \quad  for the degree $\ell=5$;

$n=(4d^3-6d^2+14d+3)/3$ \quad  for the degree $\ell =7$.

These results can be obtained with Smolyak's method, we explain 
the case $\ell=7$: 
First we take, as in \eqref{ass2}, the values 
$n_2=n_3=3$ and $n_4=7$. 
Now observe that the 4 new points of $X^4_j$ are symmetric 
but otherwise arbitrary. Hence we can take (together with 0) 
the 5-point Gau\ss\  rule with degree 9. 
This means that $2d$ weights disappear and hence $n$ is decreased 
by $2d$ compared to the general situation of \eqref{ass2}.

The best results (so far) for $\ell >7$ can be described 
in the following way: 
We use again the sequence 
$
m_i \ge 2i-1 
$ 
and so called 
``delayed Kronrod-Patterson-formulas''. 
The $n_i$ are defined as follows: 
$n_1=1$, 
$n_2=n_3=3$, 
$n_4=n_5=n_6=7$, 
$n_7= \dots  = n_{12} = 15$,
$n_{13} = \dots = n_{24} = 31$ 
and so on.
Some of these numbers are larger than $2i-1$ and hence 
we can modify those $n_i$, used by 
Petras (2003), to
$$
\tilde n_i := \min( n_i , 2i-1) .
$$
In this way one obtains the values from Table~3, 
see Genz (1986) who obtained the same results. 

\begin{table}[h]
\caption{Known values for the Lebesgue measure}
\label{ta9}
\[
\begin {array}{rrrrrr} 
\ell  & N(\ell ,5) & N (\ell ,10) & N (\ell ,15) & N (\ell ,20) 
& N (\ell ,25)  \\[\medskipamount]
3&10&20&30&40&50
\\\noalign{\medskip}
5&51&201&451&801&1\,251
\\\noalign{\medskip}
7&141&1\,181&4\,121&9\,961&19\,701
\\\noalign{\medskip}
9&391&5\,281&27\,671&90\,561&226\,951
\\\noalign{\medskip}
11&903&19\,105&145\,607&641\,409&2\,037\,511
\\\noalign{\medskip}
13&1\,733&60\,205&642\,417&3\,745\,369&14\,996\,061
\\\noalign{\medskip}
15&3\,263&168\,825&2\,473\,287&18\,743\,249&93\,755\,311
\\\noalign{\medskip}
17&5\,983&431\,265&8\,522\,247&82\,703\,329&511\,676\,911
\end {array}
\]
\end{table}

\begin{rem} 
Observe that, up to now, there is nothing better known 
than the fully symmetric formulas that were introduced 
more than 40 years ago. 
We do not claim that the results of Table~3 are optimal 
for fully symmetric (or Smolyak) rules.
It was proved by Petras (2003), however, that only minor 
improvements are possible if one uses Smolyak formulas. 
The same also holds for the more general fully symmetric 
formulas. For fixed $\ell=2k+1$ and large $d$, the number of points 
is (at least)  of the order
\begin{equation}   \label{eq35} 
N(2k+1, d, \rho)   \approx   \frac{2^k \, d^k}{k!} , 
\end{equation} 
while the lower bound of M\"oller is only of the order
$\frac{2d^k}{k!}$. 
Observe that \eqref{eq35} holds for all the versions of Smolyak's 
method that we presented here. 
\end{rem} 

\begin{rem} 
By Lemma~\ref{t2} we have the bound
$$ 
N(2k+1,d, 1 ) \leq \binom{k+d}{d} \cdot \min (2^k,2^d)
$$ 
for the Smolyak methods described here. 
For fixed $k$ also Kuperberg (2004) 
obtains a bound of the form
$$
N(2k+1,d, 1 ) \leq \binom{k+d}{d} \cdot C_k .
$$
The constant $C_k$ is of the order
$2 \cdot k^k \cdot k!$, much bigger than 
$2^k$. 
However, Kuperberg (2004) obtains cubature formulas with positive 
(even equal) weights. This is a great advantage, in
particular if the function values $f(\bx_i)$ are given only 
approximately. 

For a cubature formula $Q_n$ we define its condition number 
$$
\sigma (Q_n) = \frac{\Vert Q_n \Vert_\infty}{\Vert I^\rho_d\Vert_\infty}
= \frac{\sum_{i=1}^n |a_i|}{\int_\Omega \rho(\bx) d\bx}. 
$$
A cubature formula with positive weights has condition number 
$\sigma (Q_n)=1$ if it is exact for the constant functions. 
The known Smolyak formulas of degree 5 and 7 have a condition number of 
roughly $d^2$ and $d^3$, respectively. 
See Remark 8 which also shows that our new formulas have roughly the 
same condition numbers. 
\end{rem} 

\section{Cubature formulas for the sphere and for $M_{d,k}$}  

In the following we need some known results for cubature 
formulas for the sphere. We use these results and the Smolyak 
method to construct efficient cubature formulas for 
the linear functional
$$
M_{d,k} (f) = \sum_{ \bx \in F(d,k)} f(\bx) 
$$
where 
$$
\bx \in F(d,k) \quad  \iff \quad  x_i \in \{ \pm 1, \, 0 \}, \quad  
\sum x_i^2 = k .
$$
Of course $M_{d,k}$ itself is a cubature formula using
$2^k \binom{d}{k}$ function values, where $k \le d$. 
The point is to find a cubature formula for $M_{d,k}$ 
that is exact for polynomials from 
$\PP (2k+1,d)$ and uses only about $2 \binom{d}{k} \approx 
2d^k/k!$ points, which is the order of the lower bound of M\"oller. 


To achieve this we use two cubature formulas for the sphere
that are exact for polynomials in $\PP (2k+1,d)$.
The first formula is obtained from the Smolyak method for the
Gaussian weight function (\ref{eq16}) by projection onto the sphere
of radius $\sqrt{k}$. It has the form
\begin{equation}
\label{mdk}
 w\, M_{d,k}(f) + Q_r(f) 
\end{equation}
where $Q_r(f)$ is a cubature formula with $r=O(d^{k-1})$ points
and $w>0$. 
In particular, we can take $r \le n^{*}(d+k,d) - 2^k \binom{d}{k}$ with
$n^*(d+k,d)$ from Remark \ref{rem4}. This leads to
$$
 r \le 2d \quad \hbox{for}\quad k=2  \qquad \hbox{and} \qquad  
 r \le 2d^2 \quad \hbox{for}\quad k=3.
$$
This works for any degree $2k+1$ of exactness. 
The second formula $\widetilde{Q}_n(f)$ for $k=2,3$ is taken 
from Mysovskikh (1968),
see also Mysovskikh (1981).
It uses
$$
n=d^2+3d+2 \quad \hbox{points if}\quad k=2 \quad \hbox{and} \quad d\ge 4
$$
and
$$
n=(d^3+9d^2+14d+6)/3 \quad \hbox{points if}
\quad k=3 \quad \hbox{and} \quad d\ge 6.
$$
It follows that the formula $w^{-1}(\widetilde{Q}_n(f) - Q_r(f))$ 
is a cubature formula for
$M_{d,k}$ exact for polynomials from $\PP (2k+1,d)$ which uses at most 
\begin{equation}
 \label{skf}
  d^2 + 5d + 2 \qquad \hbox{and} \qquad  
 (d^3+15d^2+14d+6)/3 
\end{equation}     
points for $k=2$ and $k=3$, respectively.

Let us finally explain how a Smolyak formula for the Gaussian
weight function leads via projection onto the sphere $R \, \SP^{d-1}$ 
of radius $R=\sqrt{k}$ to a cubature formula of the same degree of exactness.
To this end, for $r>0$, let $\omega_r$ be the surface measure
on the sphere of radius $r$. 
Let also $P$ be the radial projection from $\R^d \setminus \{0\}$ 
onto $R \,\SP^{d-1}$ given by $P\bx = R\bx/{\Vert \bx \Vert_2^2}$. 
Furthermore, let 
\[ Q_n(f) = \sum_{i=1}^n a_i f(\bx_i) \] 
be an arbitrary  cubature formula  
which is centrally symmetric. 
Obviously, any Smolyak formula considered above has this 
property.
We assume that $Q_n$ has degree of exactness $2k+1$ for
the Gaussian weight function. 
Let $\bx^\alpha = x_1^{\alpha_1} \ldots x_d^{\alpha_d}$ be a
monomial of degree $|\alpha|=\alpha_1+\ldots+\alpha_d=2k$. 
Using polar coordinates, we obtain
\begin{eqnarray*}
 \int_{\R^d} \bx^\alpha \exp({-\Vert \bx \Vert_2^2}) d\bx
  &=& \int_0^\infty \int_{r\,\SP^{d-1}} \bx^\alpha  
  d\omega_r(\bx) e^{-r^2} dr \\
  &=& \int_0^\infty  (r/R)^{d-1+2k} e^{-r^2} dr 
  \int_{R\,\SP^{d-1}} \bx^\alpha  d\omega_R(\bx)\\
  &=& c(R,d,k) \int_{R\,\SP^{d-1}} \bx^\alpha  d\omega_R(\bx). 
\end{eqnarray*}
We also have
\[
 Q_n(\bx^\alpha) = \sum_{i=1}^n a_i \bx_i^\alpha 
 = \sum_{i=1}^n a_i (\|\bx_i\|_2/R)^{2k} (P\bx_i)^\alpha. 
\]
Whenever one of the points $\bx_i=0$, we simply drop the corresponding term. 
Since 
\[ Q_n(\bx^\alpha) = \int_{\R^d} \bx^\alpha \exp({-\Vert \bx \Vert_2^2}) d\bx,\]
we obtain that 
\[ PQ_n(\bx^\alpha) = \int_{R\,\SP^{d-1}} \bx^\alpha  d\omega_R(\bx)\]
where 
\[ PQ_n(f) = \sum_{i=1}^n b_i f(P\bx_i) \]
with
\[ b_i = \frac{a_i \|\bx_i\|_2^{2k}}{R^{2k} c(R,d,k)}. \]
So $PQ_n(f)$ is a cubature formula for the sphere 
$R \, \SP^{d-1}$  which is exact for
homogeneous polynomials of degree $2k$. 
Since it inherits the central symmetry from $Q_n$,
it is also exact for homogeneous polynomials of degree $2k+1$. 
Since any polynomial in
$\PP (2k+1,d)$ restricted to $R \, \SP^{d-1}$ 
is a sum of two homogeneous polynomials of degree $2k$
and $2k+1$, respectively, $PQ_n$ is exact for all such polynomials.

If we choose the sets $X_j^2$ in the construction 
of the Smolyak formula for the
Gaussian measure equal, say $X_j^2=X^2=\{-a,0,a\}$, 
then the points $\bx \in a F(d,k)$
are present in the Smolyak formula and get equal positive weights. 
So the projection of this formula to the sphere $R \, \SP^{d-1}$ 
has indeed the form (\ref{mdk}).
 
\begin{rem}
\label{rem5}
It will be important later on that the cubature formula derived for $M_{d,k}$
uses only points on the same sphere of radius $R=\sqrt{k}$ where the points 
in $F(d,k)$ live.
\end{rem} 

\section{Cubature formulas for general weight functions}

We now derive our main result which is formulated in the
following theorem.

\begin{theorem}
 \label{gen2}
Let $\Omega$ and $\varrho$ be as always and let $k=2,3$. 
In the case $k=2$ we assume $d \ge 4$, in the case $k=3$ 
we assume $d \ge 6$. 
Then there exists a cubature formula $Q_n$ for $I_d^\varrho$
 with degree $2k+1$ of exactness which uses at most
\begin{equation}   \label{neuallg} 
d^2 + 9d + 1 \qquad \hbox{and} \qquad  (d^3+33d^2+14d+3)/3
\end{equation}  
 points for $k=2$ and $k=3$, respectively. 
If the one-dimensional weight functions $\varrho_i$ are equal
(the fully symmetric case) 
then the number of points can be reduced to
\begin{equation} \label{neuspec} 
d^2 + 7d + 1 \qquad \hbox{and} \qquad  (d^3+21d^2+20d+3)/3 
\end{equation}  
for $k=2$ and $k=3$, respectively. \end{theorem}  

\begin{proof}[\bf Proof]
We start by describing how one can pass from the special 
cubature formulas for $M_{d,k}$
constructed in the preceding section to cubature formulas for
general weight functions $\varrho$ as in the introduction. 
By proper scaling, we may assume that the radius of the domain $\Omega$ of 
integration is at least $\sqrt{k}$.
First, choose a Smolyak formula $Q_m$ for 
$\varrho$ that is exact for polynomials
from $\PP (2k+1,d)$ and satisfies 
$X_1^2 = \ldots = X_d^2=\{-1,0,1\}$. 
Then $Q_m$ has the form
\begin{equation}
 \label{gen1}
  Q_m = v M_{d,k} + Q_s
\end{equation}
for some $v>0$ and
\[ s = n(k+d,d) - 2^k \binom{d}{k}.\]
In general, we have to use the case where $n_i = 2i-1$ for all $i\ge 1$.
Then we obtain 
$$
s= 4d+1  \qquad \hbox{and} \qquad s= (18d^2 +3)/3 
$$
for $k=2$ and $k=3$, respectively.
Now we replace the part $M_{d,k}$ in (\ref{gen1}) with the formula
derived in the preceding section which uses at most as much points
as given in (\ref{skf}). By Remark \ref{rem5} all points 
of the final cubature formula 
$$
\frac{v}{w} ( \wt Q_n - Q_r ) + Q_s 
$$
are in the interior of $\Omega$.
This cubature formula needs at most $n+r+s$ function values. 
This leads to cubature formulas with
\[ d^2 + 9d + 3 \qquad \hbox{and} \qquad  (d^3+33d^2+14d+9)/3 \]  
points for $k=2$ and $k=3$, respectively, which exceeds
(\ref{neuallg}) by just two knots. 

A further reduction is possible if knots of 
$\wt Q_n$, $Q_r$ and/or $Q_s$ coincide.
We explain how this leads to the reduced number of knots
in (\ref{neuspec}) in the fully symmetric case.
The reduction by two knots in the general
case is achieved similar (and easier).

To simplify notation, we denote by $M^r_{d,k}$ for $r>0$  the
cubature formula 
$$
M^r_{d,k} (f) = \sum_{ \bx \in F^r(d,k)} f(\bx) 
$$
where 
$$
\bx \in F^r(d,k) \quad  \iff \quad  x_i \in \{ \pm r, \, 0 \}, \quad  
\sum x_i^2 = k r^2 .
$$
Observe that $M^1_{d,k}=M_{d,k}$.

We further need some notation for fomulas derived from the 
simplex. Let $S$ be a regular simplex with vertices in the
unit sphere $\SP^{d-1}$. Let $S^r_{d,k}$ be the
cubature formula 
$$
S^r_{d,k} (f) = \sum_{ \bx \in  G^r(d,k)}(  f(\bx)+f(-\bx) ) 
$$
where $G^r(d,k)$ is the set of all projections of the centers
of the $(k-1)$-dimensional faces of $S$ onto the sphere of radius $r$.
For the formulas of degree 7 we need one more cubature formula. 
Denote by $p_{ij}$ the $(d+1)d$ points of the form 
$$
p_{ij} = \frac{1}{4} v_i + \frac{3}{4} v_j , 
$$
where $v_i$ and $v_j$ are different vertices of the simplex. 
Then let $H^r(d)$ be the set of all
$rp_{ij}/\Vert p_{ij} \Vert$ and define the cubature formula 
$\wt S^r_d$ by
$$
\wt S^r_d (f) = \sum_{ \bx \in  H^r(d)}(  f(\bx)+f(-\bx) ) . 
$$
Finally, let $\omega_d$ be the surface area of $\SP^{d-1}$. 

So assume now that $\varrho_1=\ldots = \varrho_d$.
We further assume without loss of generality that 
$\Omega \supset [-1,1]^d$.
We treat the degree five and seven cases separately. 

\medskip 

\noindent 
{\bf Degree five. } \ 
The projected Smolyak formula with degree of
exactness 5 for the sphere
$\SP^{d-1}$ with $d\ge 3$ needs $2d^2$ points and has the form
\begin{equation}  \label{eqm} 
  u_1 \, M^1_{d,1} + u_2 \, M^{1/\sqrt{2}}_{d,2} 
\end{equation}  
with 
\[ u_1= \frac{4-d}{2 d (d+2)} \omega_d \qquad \hbox{and} \qquad  
   u_2= \frac{1}{d (d+2)} \omega_d.\]
This formula can be found  in
Stroud (1971) or as formula 11) for the sphere in Mysovskikh (1981).

The second formula with degree of
exactness 5 for the sphere
$\SP^{d-1}$ with $d\ge 4$ needs $(d+1)(d+2)$ points and has the form
\begin{equation} \label{eqs}   
  v_1 \, S^1_{d,1} + v_2 \, S^1_{d,2} 
\end{equation}  
with 
\[ v_1= \frac{d(7-d)}{2 (d+1)^2 (d+2)} \omega_d \qquad \hbox{and} \qquad  
   v_2= \frac{2 (d-1)^2}{d (d+1)^2 (d+2)} \omega_d.\]
This formula can be found in Mysovskikh (1968) or as formula 
7) for the sphere in Mysovskikh (1981).

Putting (\ref{eqm}) and (\ref{eqs}) together 
gives the following formula with degree of
exactness 5 for $M^{1/\sqrt{2}}_{d,2}$:
\begin{equation}
 \label{eqm2}
 \frac{1}{u_2}
 (v_1 S^{1}_{d,1} + v_2 \, S^{1}_{d,2} - u_1 M^{1}_{d,1}).
\end{equation}
    
We also need a Smolyak type formula for the weight function $\varrho$
with degree of exactness 5 which has the form
\begin{equation}
 \label{eqc}
a_1 \, M^{1/\sqrt{2}}_{d,2} + a_2 \, M^{1/\sqrt{2}}_{d,1} + a_3 \, M^{\gamma}_{d,1} 
  + a_4 Q_0,
\end{equation}
where $Q_0(f)=f(0)$ and $\gamma \in (0,1) \setminus \{ 1/\sqrt{2} \}$.
The coefficients $a_1,\ldots,a_4$ can be 
derived either from the Smolyak construction
or from direct computation using Sobolev's theorem which tells us that
our formula has the required degree of exactness if it integrates the
polynomials $1, x_1^2, x_1^4, x_1^2 x_2^2$ correctly. This leads to a linear
system of 4 equations for $a_1,\ldots,a_4$ which has a unique solution.
To minimize the number of knots we choose $\gamma=1$.
   
Finally, we replace $M^{1/\sqrt{2}}_{d,2}$ in formula (\ref{eqc}) with the
expression (\ref{eqm2}). This leads to a formula 
\begin{equation}
\label{eqfinal}
\a_1 \left(S^{1/\sqrt{2}}_{d,2} + \frac{d^2(7-d)}{4(d-1)^2}
\, S^{1/\sqrt{2}}_{d,1} \right) 
+ \a_2 \, M^{1/\sqrt{2}}_{d,1}  +  \a_3 \, M^{1/2}_{d,1} + \a_4 Q_0 
\end{equation}
which is exact of degree 5 
for integration with respect to $\varrho$ with $d\ge 4$. 
The coefficients $\alpha_1,\ldots,\alpha_4$ can be directly derived
using the polynomials $1, x_1^2, x_1^4, x_1^2 x_2^2$.
Alternatively, they are related to $a_1,\ldots,a_4$ via
\[ \alpha_1 = \frac{2 (d-1)^2}{(d+1)^2} a_1, \quad
   \alpha_2 = a_3 - \frac{4-d}{2} a_1, \quad
   \alpha_3 = a_2, \quad
   \alpha_4 = a_4.
\]
Observe that
we have chosen our formulas so that the final number of knots is
$d^2+7d+3$.
This can be further reduced to
\[ 
d^2+7d+1
\]
if we choose one of the vertices of the regular simplex $S$ as the unit vector
$(1,0,\ldots,0)$.
Observe also that in the case $d=7$ the number of knots reduces even further.

\medskip 

\noindent 
{\bf Degree seven. } \ 
Let us now derive a formula with degree of exactness 7, 
i.e.,  $k=3$. The projected Smolyak formula with degree of
exactness 7 for the sphere
$\SP^{d-1}$ with $d\ge 3$ needs $(4d^3-6d^2+8d)/3$ points and has the form
\begin{equation}
\label{eqm7}
u_1 \, M^1_{d,1} + u_2 \, M^{1/\sqrt{2}}_{d,2} + u_3 \,
M^{1/\sqrt{3}}_{d,3} . 
\end{equation}
This formula can be found  in
Stroud (1971) or as formula 21) for the sphere in Mysovskikh (1981).

The second formula with degree of
exactness 7 for the sphere
$\SP^{d-1}$ with $d\ge 6$ needs $(d^3+9d^2+14d+6)/3$ points and has the form
\begin{equation}
\label{eqs7}
v_1 \, S^1_{d,1} + v_2 \, S^1_{d,2}  + v_3 \, S^1_{d,3} + 
v_4 \, \wt S^1_d . 
\end{equation}
This formula can be found in Mysovskikh (1968) or as formula 
13) for the sphere in Mysovskikh (1981).

Putting (\ref{eqm7}) and (\ref{eqs7}) together
gives the following formula with degree of
exactness 7 for $M^{1/\sqrt{3}}_{d,3}$:
\begin{equation}
\label{eqm27}
\frac{1}{u_3} \, \left( 
v_1 \, S^1_{d,1} + v_2 \, S^1_{d,2}  + v_3 \, S^1_{d,3} + 
v_4 \, \wt S^1_d  -
u_1 \, M^1_{d,1} - u_2 \, M^{1/\sqrt{2}}_{d,2} \right) . 
\end{equation}
    
We also need a Smolyak type formula for the weight 
function $\varrho$ with degree of exactness 7
which has the form
\begin{equation}
\label{eqc7}
a_1 \, M^{1/\sqrt{3}}_{d,3} + a_2 \, M^{1/\sqrt{3}}_{d,2} + 
a_3 \, M^{1/\sqrt{3}}_{d,1} + 
a_4 \, M^{\g_1}_{d,2} + a_5 \, M^{\g_1}_{d,1} + 
a_6 \, M^{\g_2}_{d,1} + a_7 Q_0 ,
\end{equation}
where $Q_0(f)=f(0)$ and the numbers $\g_1$ and $\g_2$ and $1/\sqrt{3}$
are pairwise different, between 0 and 1. 
To minimize the number of knots in the following
we choose $\g_1=1/\sqrt{2}$ and $\g_2=1$.

Finally, we replace $M^{1/\sqrt{3}}_{d,3}$ in formula (\ref{eqc7}) with the
expression (\ref{eqm27}). This leads to a formula
of the form 
\begin{eqnarray*}
\a_1 \, \left( 
v_1 \, S^1_{d,1} + v_2 \, S^1_{d,2}  + v_3 \, S^1_{d,3} + v_4 \, \wt S^1_d 
\right) +  & \\
\nonumber
\a_2 \, M_{d,1}^1  + 
\a_3 \, M_{d,2}^{1/\sqrt{2}} + 
\a_4 \, M_{d,1}^{1/\sqrt{2}} + 
\a_5 \, M_{d,2}^{1/\sqrt{3}} + 
\a_6 \, M_{d,1}^{1/\sqrt{3}} + 
\a_7 \, Q_0 .  & 
\end{eqnarray*}
The constants $\a_1, \dots , \a_7$ can be determined
by using the 7  polynomials
$1$, $x_1^2$, $x_1^4$, $x_1^2 x_2^2$, 
$x_1^2 x_2^2 x_3^2$, $x_1^4 x_2^2$ 
and
$x_1^6$. 
Observe that
we have chosen our formulas so that the number of knots is
$$
(d^3 + 21d^2 +20 d +9 )/3 . 
$$
This can be further reduced to
$$
(d^3 + 21d^2 +20 d +3 )/3 
$$
if we choose one of the vertices of the regular simplex $S$ as the unit vector
$(1,0,\ldots,0)$.

\end{proof}

Table \ref{ta99} contains the number of function values 
for fully symmetric weight functions.
Observe that for $\ell=7$ we have to assume $d \ge 6$. 

\begin{table}[h]
\caption{New values for fully symmetric weight functions}
\label{ta99}
\[
\begin {array}{rrrrrrrr} 
\ell  & N(\ell ,5) & N (\ell ,10) & N (\ell ,15) & N (\ell ,20) 
& N (\ell ,25)  & N(\ell, 50) & N(\ell, 100)    \\[\medskipamount]
5&61&171&331&541&801&2\,851&10\,701  
\\\noalign{\medskip}
7& - &1\,101&2\,801&5\,601&9\,751&59\,501&404\,001      
\end {array}
\]
\end{table}

It is interesting to compare these values with the lower 
bound~\eqref{moeller} of M\"oller, see Table~\ref{ta98}.

\begin{table}[h]  
\caption{M\"oller's lower bound}
\label{ta98}
\[
\begin {array}{rrrrrrrr} 
\ell  & N(\ell ,5) & N (\ell ,10) & N (\ell ,15) & N (\ell ,20) 
& N (\ell ,25)  & N(\ell, 50) & N(\ell, 100)    \\[\medskipamount]
5&31&111&241&421&651&2\,551&10\,101  
\\\noalign{\medskip}
7&80&460&1\,390&3\,120&5\,900&44\,300&343\,600   
\end {array}
\]
\end{table}

\begin{rem}
For the cube $[-1,1]^d$ with Lebesgue measure, 
Tables \ref{ta05} and \ref{ta07} contain the
coefficients $a_i$ and $\alpha_i$ in the cubature formulas
(\ref{eqc}), (\ref{eqfinal}), (\ref{eqc7}). The values of $v_1,\ldots,v_4$ and 
$u_1,u_2,u_3$ for the degree 7 formula can be found in Mysovskikh (1981).
\begin{table}[h]
\caption{Coefficients for the degree 5 formulas (\ref{eqc}) and (\ref{eqfinal})}
\label{ta05}
\begin{center}
\begin{tabular}{ccccc}
$i$          & 1 & 2 & 3 & 4 \\
$2^{-d} a_i$ & \tabmath{\frac{1}{9}}  & \tabmath{\frac{22}{45} - \frac{2d}{9}} & \tabmath{\frac{1}{30}} & 
\tabmath{\frac{2d^2}{9} - \frac{37 d}{45} +1}\\
$2^{-d} \alpha_i$ \  & \  \tabmath{\frac{2 (d-1)^2}{9 (d+1)^2}} \  &
\  \tabmath{\frac{d}{18}-\frac{17}{90}} \  &  \ 
\tabmath{\frac{22}{45} - \frac{2d}{9}} \  & \  \tabmath{\frac{2d^2}{9} - \frac{37 d}{45} +1}\\
\end{tabular}
\end{center}
\end{table}
\begin{table}[h]
\caption{Coefficients for the degree 7 formula (\ref{eqc7})}
\label{ta07}
\begin{center}
\begin{tabular}{ccccc}
$i$          & 1 & 2 & 3 & 4  \\
$2^{-d} a_i$ & \tabmath{\frac{1}{8}}  & \tabmath{\frac{7}{20} - \frac{d}{4}} & 
\tabmath{\frac{23}{70} - \frac{9}{20} + \frac{d^2}{4}}  
& \tabmath{\frac{8}{45}} \\
&&&&\\
\end{tabular}
\begin{tabular}{cccc}
$i$  & 5 & 6 & 7 \\
$2^{-d} a_i$  &
\tabmath{\frac{32}{63} - \frac{16d}{45}} & \tabmath{\frac{1}{21}}& 
\tabmath{ - \frac{d^3}{6} + \frac{5d^2}{9} - \frac{659 d}{630} +1}\\
\end{tabular}
\end{center}
\end{table}
\end{rem}

\begin{rem} 
Victoir (2004) and Kuperberg (2004) describe, in particular, 
methods for $\ell =5$ and positive weights.
For $d=100$ Victoir has $n=4^{12}=16\,777\,216$ 
and this was further improved by Kuperberg 
to $n= 65\,536$ points with positive weights. 
See the discussion in Kuperberg (2004). 

For general weights the old record was $20\,001$, see
\eqref{old}.  
Our method needs $10\,701$ 
function values, the lower bound of M\"oller is $10\,101$. 
\end{rem} 

\section{Independence of the weight function}

We now use the Smolyak formulas to show that, for any fixed $k$,
the minimal number of knots needed by a cubature formula of degree $2k+1$
does not essentially depend on the weight function. Since the M\"oller
lower bound is of order $d^k$, the following theorem shows that
the difference can only be in the lower order terms.

\begin{theorem}
 \label{gen3}
 Let $\Omega^{(j)}$ and $\varrho^{(j)}$, $j=1,2$, be two regions  and 
 weight functions in $\R^d$ as described in the introduction. 
 For $k=2,3,\ldots$, define 
 \[ c_k= \frac{2^{2k}}{(k-1)!}.\] 
 Then
\[ | N_\min (2k+1,d,\varrho^{(1)}) 
- N_\min (2k+1,d,\varrho^{(2)})| \le c_k d^{k-1} \]
 for all $d\ge k$.
\end{theorem}

\begin{proof}
Without loss of generality, 
we assume that
the cube $[-1,1]^d$ is contained
in the interior of 
$\Omega^{(1)}$ and $\Omega^{(2)}$.
We choose a cubature formula $Q_n$ for $\varrho^{(1)}$ exact for polynomials
in $\PP (2k+1,d)$ with $n = N_\min (2k+1,d,\varrho^{(1)})$.
By proper scaling if necessary we may now assume that the 
knots of $Q_n$ are in the interior of $\Omega^{(2)}$.
We also choose, for $j=1,2$, Smolyak formulas 
\[ Q_{m_j}^{Smol} = w_j M_{d,k} + Q_{r_j} \]
for $\varrho^{(j)}$ of degree $2k+1$ with $w_j>0$. 
To assure their existence,
we have to work with the case $n_i=2i-1$ for all $i$. 
In this case we can also arrange that the knots of  $Q_{r_j}$
are contained in $[-1,1]^d$. 
Then, for $d\ge k$, the estimate
\begin{equation}
 \label{gen5}
 r_j \le 2^k \binom{d+k}{k} - 2^k \binom{d}{k}
\end{equation}
follows from (\ref{g24}).
Now 
\[ \frac{w_2}{w_1} ( Q_n - Q_{r_1} ) + Q_{r_2} \]
defines a cubature rule for $\varrho^{(2)}$ exact for polynomials
in $\PP (2k+1,d)$ with at most $n+r_1+r_2$ knots.
Observe that all the knots used are in the interior of $\Omega^{(2)}$.
By (\ref{gen5}), to prove the theorem it is enough to verify the
elementary inequality
\[ 2^k \binom{d+k}{k} - 2^k \binom{d}{k} \le  \frac{2^{2k}}{(k-1)!} d^{k-1}\]
for $d\ge k$, which is equivalent to
\[
 (d+k)(d+k-1)\ldots(d+1) - d (d-1) \ldots (d-k+1) \le k \, 2^k \, d^{k-1}.
\]
Since the left-hand side of this inequality does not exceed $(d+k)^k-(d-k)^k$,
this is an immediate consequence of
\[
 (d+k)^k-(d-k)^k = 2 \sum_{\stackrel{0\le i \le k}{i \ odd}} \binom{k}{i} d^{k-i} k^{i} \le 
    2 d^{k-1} k  \sum_{\stackrel{0\le i \le k}{i \ odd}} \binom{k}{i}  
       = k \, 2^k \,   d^{k-1}.
\]
\end{proof}

\begin{rem}   \label{r10} 
Similarly, it can be shown that 
\[ | N_\min (2k+1,d,\mu_d) 
- N_\min (2k+1,d,\varrho)| \le c_k d^{k-1}, \]
where $\mu_d$ is the surface measure on the sphere $\SP^{d-1}$ and
$\varrho$ is a weight function as in Theorem \ref{gen3}.
%
%
%
%
%
%
%
%
\end{rem}

\noindent
{\bf Acknowledgment.} \
We thank two anonymous referees for helpful comments.

\def\AM{Aequationes Math.\ }
\def\BAMS{Bull. Amer. Math. Soc.\ }
\def\BIT{BIT\ }
\def\C{Computing\ }
\def\CA{Constr. Approx.\ }
\def\CAD{Comput. Aided Design\ }
\def\CAGD{Computer-Aided Geom. Design\ }
\def\CJ{Computer J.\ }
\def\ICGA{IEEE Comp. Graph. Appl.\ }
\def\IMAN{Inst. Math. Applics. Numer. Anal.\ }
\def\IND{Ind. Univ. J. Math.\ }
\def\JACM{J. ACM\ }
\def\JAM{J. Analyses Math.\ }
\def\JAT{J. Approx. Th.\ }
\def\JC{J. Complexity\ } 
\def\JIMA{J. Inst. Math. Applics.\ }
\def\JMA{SIAM J. Math. Anal.\ }
\def\JMAA{J. Math. Anal. Appl.\ }
\def\JMM{J. Math. Mech\ }
\def\JMP{J. Math. Phys.\ }
\def\JOTA{J. Optimization Th. Appl.\ }
\def\JSSC{J. Sci. Stat. Comp.\ }
\def\LAA{Linear Alg. Appl.\ }
\def\MAA{Math. Anal. Appl.\ }
\def\MC{Math. Comp.\ }
\def\NM{Numer. Math.\ }
\def\NFAO{Numer. Func. Anal. Optim.\ }
\def\PAMS{Proc. Amer. Math. Soc.\ }
\def\PEMS{Proc. Edinburgh Math. Soc.\ }
\def\PJM{Pacific J. Math.\ }
\def\RMJ{Rocky Mt. J. Math.\ }
\def\SJNA{SIAM J. Numer. Anal.\ }
\def\SMA{SMA\ }
\def\TAMS{Trans. Amer. Math. Soc.\ }
\def\TOG{Trans. Graph.\ }
\def\TOMS{ACM Trans. Math. Software\ }
\def\USSR{USSR Comput. Maths. Math. Phys.\ } 

\section*{References}

\refx
Berens, H., Schmid, H. J., and Xu, Y. (1995):
Multivariate Gaussian cubature formulae. 
Arch. Math. {\bf 64}, 26--32  

\refx 
Bungartz, H.-J., Griebel, M.  (2004):
Sparse grids. 
Acta Numerica {\bf 13}, 147-269. 

\refx
Capstick, S., Keister, B. D. (1996):
Multidimensional quadrature algorithms at higher degree and/or 
dimension. 
J. of Computational Physics {\bf 123}, 267--273


\refx
Cools, R. (1997): Constructing cubature formulas: the 
science behind the art. 
Acta Numerica {\bf 6}, 1--54

\refx
Cools, R. (2003):
An encyclopedia of cubature formulas. 
J. Complexity {\bf 19}, 445--453 

\refx
Cools, R. and Haegemans, A. (1994):
An imbedded family of cubature formulae for 
$n$-dimensional product regions. 
J. Comput. Appl. Math. {\bf 51}, 251--262


\refx
Genz, A. C. (1986):
Fully symmetric interpolatory rules for multiple integrals.
SIAM J. Numer. Anal. {\bf 23}, 1273--1283

\refx
Genz, A. C., Keister, B. D. (1996):
Fully symmetric interpolatory rules for multiple integrals over 
infinite regions with Gaussian weight. 
J. Comput. Appl. Math. {\bf 71}, 299--309

\refx
Gerstner, T., Griebel, M. (1998):
Numerical integration using sparse grids. 
Numer. Algorithms {\bf 18}, 209--232  

\refx
Kuperberg, G.  (2004):
Numerical cubature using error-correcting codes.
Preprint, arXiv:math.NA/0402047 

\refx
Lu, J., Darmofal, D. L. (2004):
Higher-dimensional integration with 
Gaussian weight for applications in probabilistic design.
SIAM J. Sci. Comput. {\bf 26}, 613--624 

\refx 
Lyness, J. N. (1965a):
Symmetric integration rules for hypercubes I-III.
Math. Comp. {\bf 19},
260--276, 394--407, 625--637

\refx  
Lyness J. N. (1965b): 
Limits on the number of function evaluations 
required by certain high-dimensional integration rules 
of hypercubic symmetry.
Math. Comp. {\bf 19}, 638--643


\refx  
McNamee, J., Stenger, F. (1967): 
Construction of fully symmetric numerical integration formulas. 
Numer. Math. {\bf 10}, 327--344

\refx  
M\"oller,  H. M. (1979):
Lower bounds for the number of nodes in cubature formulae.
In: H\"ammerlin, G., ed., Numerische Integration, ISNM {\bf
45}, pp. 221-230. Birkh\"auser, Basel

\refx
Mysovkikh, I. P. (1968):
On the construction of cubature formulas 
with the smallest number of nodes. 
Soviet Math. Dokl. {\bf 9}, 277-280. 
[Russian original: Dokl. Akad. Nauk SSSR {\bf 178}, 
1252-1254.] 

\refx
Mysovskikh, I. P. (1981):
{\it Interpolatory Cubature Formulas.}  
Nauka, Moscow. [In Russian.] 


\refx
Novak, E., Ritter, K. (1996): 
High dimensional integration of smooth functions over cubes.
\NM {\bf 75}, 79--97

\refx
Novak, E., Ritter, K. (1999):
Simple cubature formulas with high polynomial exactness.
Constr. Approx. {\bf 15}, 499--522

\refx
Novak, E., Ritter, K., Schmitt, R., Steinbauer A. (1999): 
On a recent interpolatory method for high dimensional
integration. 
J. Comput. Appl. Math. {\bf 112}, 215--228  

\refx
Petras, K. (2003):
Smolyak cubature of given polynomial degree 
with few nodes for increasing dimension.
Numer. Math. {\bf 93}, 729--753

\refx 
Schmid,  H. J. (1983): Interpolatorische Kubaturformeln.
Dissertationes Mathematicae, CCXX

\refx
Smolyak, S. A. (1963): Quadrature and interpolation formulas 
for tensor products of certain classes of functions. 
Soviet Math. Dokl. {\bf 4}, 240-243 

\refx
Stroud, A. H. (1971):
{\it Approximate calculation of multiple integrals.}
Prentice-Hall, Englewood Cliffs, NJ 

\refx
Victoir, N. (2004):
Asymmetric cubature formulae with few points in high 
dimension for symmetric measures. 
SIAM J. Numer. Anal. {\bf 42}, 209--227 

\end{document}